\documentclass[11pt]{article}
\usepackage{amsfonts}

\usepackage{graphics}
\usepackage{indentfirst}
\usepackage{cite}
\usepackage{latexsym}
\usepackage{amsmath}
\usepackage{amssymb}
\usepackage[dvips]{epsfig}
\usepackage{amscd}

\newtheorem{theorem}{Theorem}[section]
\newtheorem{remark}{Remark}[section]

\newtheorem{lemma}[theorem]{Lemma}

\newtheorem{corollary}[theorem]{Corollary}

\newcommand{\n}{\rho}

\renewcommand{\div}{{\rm div}}

\newcommand{\na}{\nabla }

\newcommand{\bt}{\begin{theorem}}
\newcommand{\bl}{\begin{lemma}}
\newcommand{\el}{\end{lemma}}
\newcommand{\et}{\end{theorem}}

\newcommand{\la}{\label}

\newcommand{\ka}{\kappa}
\newcommand{\om}{\Omega}

\newcommand{\bn}{\begin{eqnarray}}
\newcommand{\en}{\end{eqnarray}}
\newcommand{\bnn}{\begin{eqnarray*}}
\newcommand{\enn}{\end{eqnarray*}}

\newcommand{\ben}{\begin{enumerate}}
\newcommand{\een}{\end{enumerate}}

\newcommand{\ba}{\begin{aligned}}
\newcommand{\ea}{\end{aligned}}
\newcommand{\be}{\begin{equation}}
\newcommand{\ee}{\end{equation}}

\def\p{\partial}
\def\lap{\triangle}
\def\lam{\lambda}
\def\ep{\varepsilon}
\def\th{\theta}
\def\rr{\mathbb{R}^{n}}

\makeatletter      
\@addtoreset{equation}{section}
\makeatother       

\title{ Large-time behavior  for   spherically symmetric flow of viscous polytropic  gas in an exterior unbounded domain with large initial data\thanks{This work is   partially supported  by  NNSFC   11301422.}}
\date{}
\author{Zhilei L{\small IANG}
\\[0.4 cm] {\normalsize  (e-mail:  zhilei0592@gmail.com)}
\\[1cm]{\normalsize  School of Economic Mathematics, }
\\{\normalsize Southwestern University of Finance and Economics, Chengdu 611130,  China}
}

\begin{document}
\maketitle

\begin{abstract} This paper deals with the   spherically symmetric
flow of compressible viscous and polytropic ideal fluid in  an
 unbounded domain exterior to a ball in   $\rr$ with   $n\ge2.$   We show   that the global solutions are time-asymptotically
stable, in the case of general  large initial data.    The key step is  to  obtain    the uniform   (both   $x$ and $t$) bounds  on density and 
temperature   from above  and below.
The proof is based on  the elaborate
energy estimates in the Lagrangian coordinates, and our conclusion  improves the previous results to include  $n=2$.
\end{abstract}

\section{Introduction}
  We study the asymptotic behavior of   spherically symmetric solutions to a polytropic
ideal model of a compressible viscous gas over an unbounded exterior domain $\Pi=\{\xi\in \rr: |\xi|>1\}$, where $n\ge2$ denotes  the spatial dimension.   The motion of a viscous polytropic ideal gas   can be described by the   equations in Eulerian coordinates (cf.\cite{bat})
\be\la{n1}
\begin{cases} \rho_t + \div(\rho \textbf{u}) = 0,\\
 \rho (\textbf{u}_t +  \textbf{u}\cdot \na \textbf{u}) + R\nabla (\n \th)=\mu\lap \textbf{u} + (\mu + \lam)\nabla\div \textbf{u},\\
 c_{v}\n (\theta_{t}+\textbf{u}\cdot \na \theta)+R\n \th\div \textbf{u}=\kappa \lap \theta  +\lambda(\div \textbf{u})^{2} +2\mu D:D,\quad \xi\in \Pi,\,t>0.
\end{cases}
\ee
Here, as usual, the unknown functions $\n,\th$ and $\textbf{u}=(u_{1},\cdot\cdot\cdot,u_{n})$ symbol the density, the absolute temperature and the velocity, respectively.  $R,c_{v}, \ka$ are given positive  constants; $\mu$ and $\lambda$ are the constant viscous coefficients satisfy  $\mu>0,\,\,2\mu+n\lambda>0;$ and  $D=D(\textbf{u})$ is  the deformation tensor,
\bnn D_{ij}=\frac{1}{2}\left(\p_{j}u_{i}+\p_{i}u_{j}\right)\quad {\rm and}\quad D:D=\sum_{i,j=1}^{n}D_{ij}^{2}.\enn

We shall consider the equations \eqref{n1} supplemented  with the initial and boundary conditions
\be\la{e03}\rho(\xi,0)=\rho_{0}(\xi),\, \textbf{u}(\xi,0)=\textbf{u}_{0}(\xi),\,  \th(\xi,0)=\th_{0}(\xi),\quad \xi\in \overline{\Pi},\ee
and \be\la{e4} \textbf{u}(\xi,t)|_{\xi\in \p\Pi}=0,\,\,\frac{\partial\th}{\p \nu}(\xi,t)|_{\xi\in \p\Pi}=0,\quad t\ge0,\ee with $\nu$ being the exterior normal vector.

If   $(\rho_{0}(\xi),\textbf{u}_{0}(\xi),\th_{0}(\xi))$   are assumed to be
spherically symmetric, i.e.,
\be\la{e5} \rho_{0}(\xi)=\hat{\rho}_{0}(r),\,\textbf{u}_{0}(\xi)=\frac{\xi}{r}\hat{u}_{0}(r),\,\th_{0}(\xi)
=\hat{\th}_{0}(r),\quad r=|\xi|\ge1,\ee
then the  symmetric functions   $(\hat{\rho},\hat{u},\hat{\th})(r,t)$ are the unique solution  because  eqs \eqref{n1} are rotationally invariant (cf.\cite{ita}), and thereby, the \eqref{n1} takes the form   (ignore   the  " $\hat{}$ ")
\be\la{rn1}
\ba &\rho_t +\frac{(r^{n-1} \n u)_{r}}{r^{n-1}}= 0,\\
 &\rho (u_t + u\p_{r}u) + R\p_{r}( \n \th) =\beta\left(\frac{(r^{n-1}u)_{r}}{r^{n-1}}\right)_{r},\\
& c_{v}\rho (\th_t + u\p_{r}\th) +R\n \th \frac{(r^{n-1}u)_{r}}{r^{n-1}}\\
&=\kappa \frac{(r^{n-1}\th_{r})_{r}}{r^{n-1}} +\lambda\left(\frac{(r^{n-1}u)_{r}}{r^{n-1}}\right)^{2} +2\mu(\p_{r}u)^{2}+2\mu\frac{n-1}{r^{2}}u^{2},\quad r\in (1,\infty),\,t>0,
\ea
\ee
where $\beta=2\mu+\lambda>0,$  the initial and boundary conditions \eqref{e03}-\eqref{e4} become
\be\la{e6}\ba \n(r,0)=\n_{0}(r),\,u(r,0)=u_{0}(r),\,\th(r,0)=\th_{0}(r),\quad r\ge1,\ea\ee
and \be\la{e7}\ba  u(1,t)=0,\,\p_{r}\th(1,t)=0,\quad t\ge0.\ea\ee

For our analysis  convenience, it is desirable to  convert the \eqref{rn1} from the Euler coordinates $(r,t)$ into that in Lagrangian coordinates $(x,t)$. Define
\be\la{f4} r(x,t)=r_{0}(x)+\int_{0}^{t}u(r(x,\tau),\tau)d\tau,\ee
with  \be\la{f3}\int_{1}^{r_{0}(x)}y^{n-1}\n_{0}(y)dy=x.\ee
Using \eqref{f4}, \eqref{f3}, $\eqref{rn1}_{1},$
and the boundary condition $u(1,t)=0,$ we check   for  $t\ge0$
\be\la{f5}\int_{1}^{r(x,t)}y^{n-1}\rho(y,t)dy=\int_{1}^{r_{0}(x)}y^{n-1}\n_{0}(y)dy=x.\ee
By this,  $r=1$ iff $x=0$ and $r\rightarrow\infty$ iff $x\rightarrow\infty$, as long as $\rho>0$ for all $(y,t)\in [0,\infty)\times [0,\infty).$
 Moreover, it is easy to see from \eqref{f4} and \eqref{f5} that
\be\la{f9} \p_{t}r(x,t)=u(r(x,t),t)\,\,{\rm and}\,\,r^{n-1}(x,t)\rho(r(x,t),t)\p_{x}r(x,t)=1.\ee
We introduce   \be\la{f8} \tilde{v}(x,t)=:1/\rho(r(x,t),t),\,\tilde{u}(x,t)=:u(r(x,t),t),\,\,\tilde{\th}(x,t)=:\th(r(x,t),t),\ee
and express  \eqref{rn1} in terms of $(\tilde{v},\tilde{u},\tilde{\th})$ ( denoted still by $(v,u,\th)$ below) in variables  $(x,t)$
\be\la{eq}
\ba &v_t =(r^{n-1}u)_{x},\\
 &u_t   = r^{n-1} \sigma_{x}, \\
& c_{v} \theta_{t} =\kappa \left(\frac{r^{2(n-1)}\theta_{x}}{v}\right)_{x} +  (r^{n-1} u)_{x}\sigma -2\mu(n-1)(r^{n-2}u^{2})_{x}, \quad x\in \Omega,\,t>0,
\ea
\ee
 where   $\sigma= \beta (r^{n-1}u)_{x}/v-R\th/v,\,\, \Omega=(0,+\infty),$   the initial functions
\be\la{f6} v(x,0)=v_{0}(x),\,u(x,0)=u_{0}(x),\,\th(x,0)=\th_{0}(x),\quad x\in \Omega,\ee
  the boundary and the far field behavior \be\la{f7}u(0,t)=0,\,\,\p_{x}\th(0,t)=0,\quad  \lim_{x\rightarrow\infty}(v(x,t), u(x,t),\th(x,t))=(1,0,1)\quad t\ge 0.\ee

In view of
\eqref{f8},   we reduce the \eqref{f4} and \eqref{f9} to
\be\la{2.3} r(x,t)=r_{0}(x)+\int_{0}^{t}u(x,s)ds,\quad r_{t}=u,\quad r^{n-1}r_{x}=v.\ee
Integrating the last equality in \eqref{2.3}   yields
\be\la{2.1}r^{n}(x,t)=1+n\int_{0}^{x}v(y,t)dy.\ee
Furthermore,  it follows from    \cite[eq.(3.19)]{jiang1} that
\be\la{2.8} r(x,t)\ge r(0,t)=1,\quad (x,t)\in \overline{\Omega}\times [0,\infty).\ee

A lot of works have been done on the existence, stability and large time behavior
of solutions to the compressible Navier-Stokes equations for either isentropic or
non-isentropic flow, and the progress is much satisfactory in  case of the assumption that the  initial data having a small oscillatory around a  non-vacuum equilibrium. See  \cite{hoff,ita,jiang,jiang1,jiang2,jiang3,kan,ll,lz,mats,mn,nn1,nn,ok,qin} and  cited therein.

The first work on the global existence in the 1-dimensional (1-D) viscous polytropic ideal gas for large initial data was due to Kazhikhov-Shelukhin \cite{kazhikhov}.
 For the    initial-boundary-value (IBV) problem \eqref{eq}-\eqref{f7}, the
  global   solution     has been studied for a bounded annular domain (cf.\cite{nik,yashi,yashi1}).
By considering   an approximate problem in bounded domain and obtaining the   a priori estimates independent of annular domain,  Jiang \cite{jiang1} established the unique  global solution  to the IBV  problem \eqref{eq}-\eqref{f7} in unbounded  exterior domain.

\begin{theorem}[Global existence in \cite{jiang1}] \la{t1}  Assume that the  initial function  in \eqref{f6} satisfy
\be\ba\la{e2}  v_{0}-1,\,u_{0},\,\th_{0}-1,\,r^{n-1}\p_{x}v_{0},\,r^{n-1}\p_{x}u_{0},\,r^{n-1}\p_{x}\th_{0} \in L^{2}(\Omega), \ea\ee
\be\ba\la{e3}  \inf_{x\in \overline{\Omega}}v_{0}(x)>0\quad {\rm and}\quad\inf_{x\in\overline{\Omega}}\th_{0}(x)>0,
 \ea\ee  and  are compatible with the boundary conditions \eqref{f7}.   Then for any fixed $T>0,$  the problem  \eqref{eq}-\eqref{f7} admits   a unique global (large)
generalized solution $(v,u,\th)$   over $[0,T]$, with $v(x, t)$ and $\th(x, t)$ having positive bounds from above and below (depending on $T$). Moreover,
\bnn\ba  &v-1,\,u,\,\th-1 \in L^{\infty}\left(0,T;H^{1}(\Omega)\right),\,\,v_{t},\, r^{n-1}u_{x},\, r^{n-1}\th_{x}\in L^{\infty}\left(0,T;L^{2}(\Omega)\right),\\
 &v_{xt},\, u_{t},\, \th_{t},\,  r^{2(n-1)} u_{xx},\, r^{2(n-1)}\th_{xx} \in L^{2}\left(0,T;L^{2}(\Omega)\right).\ea\enn
 \end{theorem}

However, the study on   asymptotics    of large solutions
to  the problem   \eqref{eq}-\eqref{f7} is less satisfactory.  In dimension $n\ge 3,$
Jiang \cite{jiang1}  obtained the bound on $v(x,t)$ independent of $t$, and obtained partial results on the asymptotic behavior of solutions,  but leaves the (upper) bound on $\th(x,t)$ open.   Later,  Nakamura-Nishibata    considered the similar  problem  (driven by  small force).  They   \cite{nn}  succeeded in obtaining  the uniform bounds on both $v$ and  $\th$, and consequently,  the corresponding  stationary solution is shown to be time-asymptotically stable.
Unfortunately,  the results  in  \cite{jiang1,nn} do not hold true  for dimension  $n=2$ due to some technical difficulties.
 
In this current paper, we discuss the   large-time   behavior  of   solutions to the IBV problem \eqref{eq}-\eqref{f7} in any multi-dimensional case (includeing $n=2$).   Of course, a key step is deriving uniform bounds for both $v$ and $\th$ in terms of  $x$ and $t$ variables.  Different from previous results in this topic,     some   new  technique  is needed for $n=2$.  Firstly,  to  deal with the bounds for $\th$,    we have the essential  observation: \emph{the spatial  domain keeps bounded once the temperature    $\th(x,t)$ is far away from the equilibrium state (i.e., $\th\equiv1)$.} In particular,   we multiply the energy equation $\eqref{eq}_{3}$ by $(\th-2)_{+}$ and obtain the critical estimates on $\int_{0}^{\infty}\|r^{n-1}\th_{x}\|_{L^{2}(\Omega)}^{2}dt.$  See Lemma \ref{l3.2} below for detailed proof.  Secondly,  we   borrow  some ideas developed  by Jiang  \cite{jiang2,jiang3} and  localize   the  classical  representation (see, e.g.,\cite{kazhi}), but some modifications  are  also needed to handle the case $n=2.$

The    theorem below states our main results.

\begin{theorem}[Large-time behavior]\la{t2} Under the same assumptions made on the initial functions,  the global   solution   $(v, u, \th)$ described in Theorem \ref{t1}  satisfies
\be\la{a2}\ba  \sup_{t\in [0,\infty)}&\left(\|(v-1, u, \th-1)(\cdot,t)\|_{L^{2}(\Omega)}+\|r^{n-1}(v_{x}, u_{x}, \th_{x})(\cdot,t)\|_{L^{2}(\Omega)}\right) \\ &\quad+\int_{0}^{\infty}\left(\|r^{2(n-1)}u_{xx}\|_{L^{2}(\Omega)}^{2}+\|r^{2(n-1)}\th_{xx}\|_{L^{2}(\Omega)}^{2}\right)dt\le C,
\ea\ee
and \be\la{a1}\ba C^{-1}\le v(x,t),\quad  \th(x,t)\le C,\quad \forall\,\,\,(x,t)\in \overline{\Omega}\times [0,\infty),
\ea\ee where the $C$ depends only on $\mu,\lambda,  R ,c_{v},\ka, n,$   and the initial data.
Moreover,   \be\la{a03}\ba  \lim_{t\rightarrow\infty}\|(v-1, u, \th-1)(\cdot,t)\|_{C(\overline{\Omega})}=0.
\ea\ee
 \end{theorem}

\begin{remark} Our method applies to the  boundary conditions
\bnn u(0,t)=0,\quad \th(0,t)=1,\quad t\ge0.\enn \end{remark}

Throughout this paper, $C(\overline{\Omega}), L^{p}(\Omega)$ and $ H^{1}(\Omega)$ denote  the usual Sobolev spaces. See, e.g., the definitions in \cite{ad}.   The    $C,C_{i}>1\quad(i=1,2,3)$ are  generic  constants which may rely on $\mu, \lambda, R, c_{v}, \ka, n,$ and the initial data,  but does not depend on the time $T$. In addition,  we use  $C(\alpha)$  to  emphasize  that $C$ depends on $\alpha.$

The rest of this paper is  arranged  as follows:
 In section 2,  some known Lemmas and facts are collected  for proving Theorem \ref{t2}. 
Sections 3-4 are devoted to deriving the bounds on  $v$ and $\th$.
In Section 5, we aim to prove some needed  $L^{2}$-norm  estimates on   derivatives,  and  in the final Section 6, we complete  the proof of  Theorem \ref{t2}.

\section{Preliminaries}
The first lemma is responsible for the basic energy estimate of the solutions, whose proof  is available in, e.g., \cite{nn1}. 
\begin{lemma} \la{lem0}   The
solution $(v,u,\th)$ obtained  in Theorem \ref{t1} satisfies
\be\ba\la{2.4} \sup_{0\le t\le \infty}\int_{\Omega} U(x,t)
 +\int_{0}^{\infty}\int_{\Omega}\left(\frac{v u^{2}}{r^{2}\th}+\frac{r^{2(n-1)}u_{x}^{2}}{v\th}+\frac{(r^{n-1}u)_{x}^{2}}{v\th} +\frac{r^{2(n-1)}\th_{x}^{2}}{v\th^{2}}\right)
 \le C,\ea\ee where \bnn U =\left(R(v-\ln v-1)+\frac{1}{2}u^{2}+c_{v}(\th-\ln \th-1)\right) .\enn\end{lemma}

With the help of \eqref{2.4},    Jensen's inequality  guarantees
\bnn \int_{k}^{k+1}v-\ln \int_{k}^{k+1}v-1,\,\,\,\int_{k}^{k+1}\th-\ln \int_{k}^{k+1}\th-1\le C,\quad k=0,1,2,\cdot\cdot\cdot,\enn
and thus
\be\la{2.5}0<\alpha_{1}\le v(a_{k}(t),t)=\int_{k}^{k+1}v(x,t),\quad \th(b_{k}(t),t)=\int_{k}^{k+1}\th(x,t)\le \alpha_{2}<\infty,\ee
where $\alpha_{1},\, \alpha_{2}$ are two   roots of the equation  $y-\ln y-1=C.$

\section{Uniform bounds of $v(x,t)$}

\begin{lemma}\la{lem1} Let   $(v,u,\th)$ be the   solution described in Theorem \ref{t1}. Then it satisfies
\be\la{e8} C^{-1}\le v(x,t)\le C,\quad (x,t)\in\overline{\Omega}\times [0,+\infty).\ee \end{lemma}
\emph{Proof.} The process is divided into several steps.

\emph{\underline{Local representation for $v(x,t)$}:}

Let
\bnn\ba \varphi(x)=\left\{
           \begin{array}{ll}
            1, & x\le k; \\
            k+1-x, & k\le x\le k+1;\\
            0, & x\ge k+1.
           \end{array}
         \right.
\ea\enn
Make use of  \eqref{2.3} and $\eqref{eq}_{1}$,  multiply $\eqref{eq}_{2}$   by $\varphi$ and after integration by parts, we receive  \be\la{3.5}\ba  v(x,t)=B(x,t)Y(x,t) +\frac{R}{\beta}\int_{0}^{t} \frac{\th(x,\tau)B(x,t)Y(x,t)}{B(x,\tau)Y(x,\tau)}d\tau,\quad x\in I_{k}, t\ge 0,\ea\ee
where  \bnn B(x,t)= v_{0} \exp\left\{\frac{1}{\beta} \int_{x}^{\infty} \varphi (r_{0}^{1-n}u_{0}-r^{1-n}u)\right\},\enn
 \bnn Y(x,t)=\exp\left\{\frac{1}{\beta}\left(\int_{0}^{t}\int_{I_{k+1}}  \sigma- (n-1)\int_{0}^{t}\int_{x}^{\infty}\varphi r^{-n}u^{2}\right)\right\}\enn with $ I_{k}=(k-1,k)$ and $k\in \mathbb{N}_{+}.$ 
 
\emph{\underline{Upper  bound of $v(x,t)$:}}

It is easy to check  from \eqref{2.8} and  \eqref{2.4} that
\be\la{3.6} C^{-1}\le B(x,t)\le C,\quad x\in I_{k}.\ee
Next to handle $Y(x,t)$.  In view of \eqref{2.8}, \eqref{2.4} and \eqref{2.5},  the same argument as    \cite[Lemma 2.4]{jiang2} shows 
\bnn\ba  -\int_{s}^{t} \inf_{x\in I_{k+1}} \th(\cdot,\tau)\le \left\{
                                                                       \begin{array}{ll}
                                                                         0, & 0\le t-s\le 1, \\
                                                                         -C(t-s), & 1\le t-s.
                                                                       \end{array}
                                                                     \right.
\ea\enn
This, along with \eqref{2.4},   \eqref{2.5} and   Jensen's inequality, yields
\bnn\ba   &\int_{s}^{t}\int_{I_{k+1}} \sigma -(n-1)\int_{0}^{t}\int_{x}^{\infty}\varphi r^{-n}u^{2}\\
 &\le \int_{s}^{t}\int_{I_{k+1}}  \sigma  \le C\int_{s}^{t}\int_{I_{k+1}} \frac{(r^{n-1}u)_{x}^{2}}{v\th}-\frac{R}{2}\int_{s}^{t}\int_{I_{k+1}} \frac{\theta}{v}\\
&\le C-\frac{R}{2}\int_{s}^{t}\inf_{I_{k+1}} \th(\cdot,\tau)\left(\int_{I_{k+1}} v\right)^{-1}\le C-C(t-s).
\ea\enn
Hence, for any $x\in I_{k}$,
\be\la{3.10} 0\le Y(x,t)/Y(x,s)\le C\exp\{-C(t-s)\},\quad 0\le s< t.\ee

Thanks  to  \eqref{3.6} and \eqref{3.10}, we deduce from  \eqref{3.5} that
\be\la{3.11}\ba  v(x,t)\le C+C \int_{0}^{t}\th(x,s)\exp\{-C(t-s)\}ds,\quad  x\in I_{k}.\ea\ee
By \eqref{2.8},  \eqref{2.5},  Holder inequality, we infer(see, e.g., \cite{jiang1})
\be\ba\la{3.12}\frac{\alpha_{1}}{2}-\alpha_{2}f(t)\max_{x\in I_{k}}  v(\cdot,t) \le \th(x,t) \le 2\alpha_{2}+2\alpha_{2}f(t)\max_{x\in I_{k}} v(\cdot,t),
\ea\ee
where \be\la{a6}f(t)=\int_{\Omega} \frac{r^{2(n-1)}\th_{x}^{2}}{v\th^{2}}.\ee
Insert  \eqref{3.12} into \eqref{3.11},  and then apply the   \eqref{2.4} and  Gronwall's inequality to find  \be\la{3.14} v(x,t)\le C,\quad (x,t)\in I_{k}\times [0,\infty),\ee for some $C$ independent of $k.$

\emph{\underline{The  bound  of $v(x,t)$ from  below:}}

 By Jensen's inequality, \eqref{3.14}, and \eqref{2.8},
\bnn \ba&\int_{I_{k+1}} \th dr-\ln \int_{I_{k+1}} \th dr-1\le \int_{I_{k+1}}\left(\th-\ln\th-1\right)dr\\
&\le C\int_{I_{k+1}}\left(\th-\ln\th-1\right)\frac{r^{n-1}}{v}dr =C\int_{I_{k+1}}\left(\th-\ln\th-1\right)dx\le C,\ea\enn
hence, \be\la{x01}C^{-1}\le \int_{I_{k+1}} \th dr\le  C.\ee

Noting   \eqref{2.3} and  \bnn\ba |r^{-n/2}u|(\cdot,t)&\le \int_{0}^{\infty}\left(|r^{-n/2}u_{x}|+|\frac{-n}{2}r^{-3n/2} vu|\right)dx, \ea\enn
we compute
\be\la{x02}\left\|\frac{u}{r^{n/2}} \right\|_{L^{\infty}}^{2} \le \left(\int_{0}^{\infty} \frac{r^{2(n-1)}u_{x}^{2}}{v\th }\right)\left(\int_{0}^{\infty} \frac{v\th}{r^{3n-2}}\right)+C\left(\int_{0}^{\infty} \frac{v u^{2}}{r^{2}\th}\right)\left(\int_{0}^{\infty} \frac{v\th}{r^{3n-2}}\right).\ee 
It follows  from \eqref{2.3} and \eqref{x01}  that
\be\ba\la{x03}  \int_{0}^{\infty} \frac{v\th}{r^{3n-2}}dx=\int_{0}^{\infty} \frac{\th}{r^{2n-1}}\frac{v}{r^{n-1}}dx=\int_{1}^{\infty} \frac{\th}{r^{2n-1}} dr=\sum_{k=2}^{\infty}\frac{1}{k^{2n-1}}\int_{I_{k}}\th dr\le C.\ea\ee 
Therefore, thanks to  \eqref{x02},\eqref{x03}, as well as \eqref{2.4},  it has \be\ba\la{x04} \left|\int_{0}^{t}\int_{x}^{\infty}\varphi r^{-n}u^{2}\right|\le C\int_{0}^{t}\|r^{-n/2}u \|_{L^{\infty}}^{2}\le C.\ea\ee

Next,  by  $\eqref{eq}_{1}$, \bnn\ba\int_{0}^{t}\int_{I_{k+1}}\frac{(r^{n-1}u)_{x}}{v}=\int_{0}^{t}\int_{I_{k+1}}\frac{v_{t}}{v}=\int_{I_{k+1}}\ln \frac{v}{v_{0}},\ea\enn
which together with  \eqref{2.5} imply 
 \bnn\ba -C_{2}\le -C_{1}-\int_{I_{k+1}}\left(v-\ln v-1\right)\le \int_{I_{k+1}}\ln \frac{v}{v_{0}}\le C_{1}+\ln \int_{I_{k+1}}v\le C_{2}.\ea\enn
The last two inequalities show that
\be\la{x05} \left|\int_{0}^{t}\int_{I_{k+1}} (r^{n-1}u)_{x}/v\right|\le C\int_{0}^{t}\|r^{-n/2}u \|_{L^{\infty}}^{2}\le C.\ee

Having \eqref{x04} and  \eqref{x05} in hand,   recalling  that $\sigma= \beta (r^{n-1}u)_{x}/v-R\th/v$, we infer 
\be\la{x06}   C^{-1}\exp\left\{\int_{s}^{t}\int_{I_{k+1}}R\th/v\right\}\le \frac{Y(x,t)}{ Y(x,s)}\le C\exp\left\{\int_{s}^{t}\int_{I_{k+1}}R\th/v\right\}.\ee
With the help of  \eqref{2.5}, \eqref{3.10}, and  \eqref{x06},  integrating  \eqref{3.5} over $I_{k}$ gives rise to
\bnn\alpha_{1}\le C\exp\{-t/C_{1}\}+C \int_{0}^{t}\exp\left\{\int_{\tau}^{t}\int_{I_{k+1}}R\th/v\right\}d\tau,\enn
which,  along  with   \eqref{3.5}-\eqref{3.10}, \eqref{3.12}, \eqref{3.14},  \eqref{x06},  deduces for $x\in I_{k}$
\be\ba\la{3.16}  v(x,t)&\ge C\int_{0}^{t} \th(x,\tau)\frac{Y(x,t)}{ Y(x,\tau)}d\tau\\
&\ge C\int_{0}^{t} \exp\left\{\int_{\tau}^{t}\int_{I_{k+1}}R\th/v\right\}d\tau-C\int_{0}^{t}f(\tau)\frac{Y(x,t)}{ Y(x,\tau)}d\tau\\
&\ge \alpha_{1}-C\exp\{-t/C_{1}\}-C\left(\int_{0}^{t/2}+\int_{t/2}^{t}\right)f(\tau)\exp\{-C(t-\tau)\}d\tau\\
&\ge \alpha_{1}-C\exp\{-t/C_{1}\}-C_{2}\exp\{-C_{3}t/2\}-\int_{t/2}^{t}f(\tau)d\tau \ge \alpha_{1}/2, \ea\ee
for all  $t\ge T_{0}$ if  $T_{0}$ large enough.
Finally, it satisfies from   \cite[eq.(4.9)]{jiang1} that
\be\la{a74} v(x,t)\ge C(T_{0}),\quad (x,t)\in \overline{\Omega}\times [0,T_{0}].\ee
The proof ends up with   \eqref{3.14}, \eqref{3.16} and \eqref{a74}.  $\Box$
\begin{corollary} Inequalities \eqref{2.1} and   \eqref{e8} show
\be\la{o4}  1+\frac{x}{C} \le r^{n}(x,t)\le  1+Cx.\ee \end{corollary}

\section{Uniform bound for $\th(x,t)$ from above}

Notice that   the set   $$\Omega_{a}(t)=\{x\in \Omega:\,\,\theta(x,t)>a>1\} $$
 is uniformly  bounded  in time, that is,    for any $t\in [0,\infty)$
\be\ba\la{z12}{\rm meas}\,\Omega_{a}(t) \le \int_{\Omega_{a}(t)}\le C(a)\int_{\Omega_{a}(t)}c_{v}(\th-\ln \th -1)\le C(a),\ea\ee by \eqref{2.4}.
This combining  with \eqref{2.5} yields
\be\la{k5}\int_{\Omega_{a}(t)}\th(x,t)\le C(a).\ee

\begin{lemma} \label{l3.2}  Let $(v,u,\th)$ be the solution described in Theorem \ref{t1}.   Then 
\be\la{e}  \sup_{0\le t\le T}\int_{\Omega}\left[(\theta-1)^{2}+u^{4}\right](x,t)+ \int_{0}^{T} \int_{\Omega}\left[(1+\th+u^{2})(r^{n-1} u)_{x}^{2}  +
       r^{2(n-1)}\th_{x}^{2}\right]
  \leq C.\ee
 \end{lemma}
\emph{Proof.}
We divide the proof   into three steps.

\emph{Step 1.}  Multiply   $\eqref{eq}_{3}$ by $(\th-2)_{+}=\max\{0,\th-2\}$, to discover
    \bnn\ba &\frac{c_{v}}{2} \int_{\Omega}(\theta-2)_{+}^{2}(x,T) +\kappa\int_{0}^{T}\int_{\Omega}\frac{r^{2(n-1)}|\p_{x}(\theta-2)_{+}|^{2}}{v}\\
& =\frac{c_{v}}{2}\int_{\Omega}(\theta_{0}-2)_{+}^{2}+2\mu(n-1)\int_{0}^{T}\int_{\Omega}r^{n-2}u^{2}\p_{x}(\theta-2)_{+}\\
&\quad+\beta\int_{0}^{T}\int_{\Omega} \frac{(r^{n-1}u)_{x}^{2}}{v}(\theta-2)_{+}-\int_{0}^{T}\int_{\Omega}\frac{R\th}{v}(r^{n-1}u)_{x}(\theta-2)_{+}. \ea\enn  Next,  multiplying  $\eqref{eq}_{2}$ by $2u(\th-2)_{+}$ gives  
 \bnn\ba  &2\int_{0}^{T}\int_{\Omega} u_{t}u(\theta-2)_{+}+2\beta\int_{0}^{T}\int_{\Omega}\frac{(r^{n-1}u)_{x}^{2}}{v}(\theta-2)_{+}\\
&=-2\beta\int_{0}^{T}\int_{\Omega}\frac{(r^{n-1}u)_{x}}{v}r^{n-1}u\p_{x}(\theta-2)_{+} -2\int_{0}^{T}\int_{\Omega}\left(\frac{R\th}{v}\right)_{x}r^{n-1}u(\theta-2)_{+}.  \ea\enn
The combination of  last two equalities arrives that
     \be\ba \label{1.3-1}
 &\frac{c_{v}}{2}\int_{\Omega}(\theta-2)_{+}^{2}+\beta\int_{0}^{T}\int_{\Omega} \frac{(r^{n-1}u)_{x}^{2}}{v}(\theta-2)_{+}+\kappa\int_{0}^{T}\int_{\Omega}\frac{r^{2(n-1)}|\p_{x}(\theta-2)_{+}|^{2}}{v}\\
     &= \frac{c_{v}}{2}\int_{\Omega}(\theta_{0}-2)_{+}^{2}+2\mu(n-1)\int_{0}^{T}\int_{\Omega} r^{n-2}u^{2}\p_{x}(\theta-2)_{+}\\
&\quad-2\beta\int_{0}^{T}\int_{\Omega}\frac{(r^{n-1}u)_{x}}{v}r^{n-1}u\p_{x}(\theta-2)_{+}
+2\int_{0}^{T}\int_{\Omega} \frac{R\th}{v}r^{n-1}u\p_{x}(\theta-2)_{+}\\
&\quad+\int_{0}^{T}\int_{\Omega}\frac{R\th}{v}(r^{n-1}u)_{x}(\theta-2)_{+}-2\int_{0}^{T}\int_{\Omega} u_{t}u(\theta-2)_{+}\\
&= \frac{c_{v}}{2}\int_{\Omega}(\theta_{0}-2)_{+}^{2}+\sum_{i=1}^{5}I_{i}. \ea\ee

The terms $I_{i}\,(i=1\sim5)$ are estimated as follows:

  By the  Cauchy-Schwarz inequality and \eqref{2.8}, the first two terms
  \bnn\ba
I_{1}+I_{2}
&\le\frac{\kappa}{4}\int_{0}^{T}\int_{\Omega}\frac{r^{2(n-1)}|\p_{x}(\theta-2)_{+}|^{2}}{v}\\
&\quad+C\int_{0}^{T}\int_{\Omega_{2}(t)} u^{4}+C\int_{0}^{T}\int_{\Omega}(r^{n-1}u)_{x}^{2}u^{2}.
\ea\enn

Second,  \be\ba\la{as2}
 I_{3}
&=2\int_{0}^{T}\int_{\Omega} \left( \frac{R(\theta-2)_{+}}{v}+ \frac{2R }{v}\right)r^{n-1}u \p_{x}(\theta-2)_{+} \\
 &=2\int_{0}^{T}\int_{\Omega}\frac{ R(\theta-2)_{+}}{v}r^{n-1}u\p_{x}(\theta-2)_{+}\\
 &\quad-4\int_{0}^{T}\int_{\Omega} \left(\frac{R}{v}\right)_{x}r^{n-1}u(\theta-2)_{+}-4\int_{0}^{T}\int_{\Omega} R\frac{(r^{n-1}u)_{x}}{v} (\theta-2)_{+}\\
&=2\int_{0}^{T}\int_{\Omega} R\frac{(\theta-2)_{+}}{v}r^{n-1}u\p_{x}(\theta-2)_{+}+I^{1}_{3}+I^{2}_{3}.\ea\ee
For one hand,  by   \eqref{e8} and integration by parts,  
 \be\ba\la{ll} I_{3}^{1}
&=4\int_{0}^{T}\int_{\Omega}  \frac{R(1-v)}{v}(r^{n-1}u)_{x}(\theta-2)_{+} +4\int_{0}^{T}\int_{\Omega} \frac{R(1-v)}{v}r^{n-1}u\p_{x}(\theta-2)_{+}\\
 &\le \beta\int_{0}^{T}\int_{\Omega}\frac{(r^{n-1}u)_{x}^{2}}{v}+C(\beta)\int_{0}^{T}\int_{\Omega}(v-1)^{2}(\theta-2)_{+}^{2}\\
 &\quad+\frac{\kappa}{16}\int_{0}^{T}\int_{\Omega}\frac{r^{2(n-1)}|\p_{x}(\theta-2)_{+}|^{2}}{v}+C(\ka)\int_{0}^{T}\int_{\Omega_{2}(t)} (v-1)^{2}u^{2}(\theta-1),\ea\ee where    $(\theta(x,t)-1)\ge 1$ in $\Omega_{2}(t).$
For another, by   $\eqref{eq}_{3}$,   \be\ba\la{ll1}
I_{3}^{2}&=-4\int_{0}^{T}\int_{\Omega} R\frac{(r^{n-1}u)_{x}}{v} (\theta-2)_{+}\\
 &=4c_{v}\int_{0}^{T}\int_{\Omega}\theta_{t}\left(1-\frac{2}{\theta}\right)_{+}+8\int_{0}^{T}\int_{\Omega_{2}(t)}  \left(\kappa\frac{r^{2(n-1)}\theta_{x}^{2}}{v\theta^{2}}+\beta\frac{(r^{n-1}u)_{x}^{2}}{v\theta}\right)\\
 &\quad -4\beta\int_{0}^{T}\int_{\Omega_{2}(t)} \frac{(r^{n-1}u)_{x}^{2}}{v}-16\mu(n-1)\int_{0}^{T}\int_{\Omega_{2}(t)}\frac{r^{n-2}u^{2}\th_{x}}{\th^{2}}\\
 &\le C+C\int_{0}^{T}\max_{x\in \Omega_{2}(t)}u^{4}-4\beta\int_{0}^{T}\int_{\Omega}\frac{(r^{n-1}u)_{x}^{2}}{v},\ea\ee
in which   the following   inequalities have been used: \bnn\ba &\int_{0}^{T}\int_{\Omega}\theta_{t}\left(1-\frac{2}{\theta}\right)_{+}\\
&=\int_{\Omega}(\theta-2\ln \theta-2(1-\ln 2))_{+}-\int_{\Omega}(\theta_{0}-2\ln \theta_{0}-2(1-\ln 2))_{+}
\le C, \ea\enn
  \bnn\ba  -\int_{0}^{T}\int_{\Omega_{2}(t)}\frac{(r^{n-1}u)_{x}^{2}}{v}
&=\int_{0}^{T}\int_{\Omega\backslash\Omega_{2}(t)}\frac{(r^{n-1}u)_{x}^{2}}{v}-\int_{0}^{T}\int_{\Omega} \frac{(r^{n-1}u)_{x}^{2}}{v}\\
&\le C-\int_{0}^{T}\int_{\Omega}\frac{(r^{n-1}u)_{x}^{2}}{v}\ea\enn
and \bnn\ba \left| \int_{0}^{T}\int_{\Omega_{2}(t)}\frac{r^{n-2}u^{2}\th_{x}}{\th^{2}}\right|&\le C\int_{0}^{T}\int_{\Omega_{2}(t)}\left(\frac{r^{2(n-1)}\th_{x}^{2}}{v\th^{2}}+\frac{u^{4}}{r^{2}\th^{2}}\right) \le C+C\int_{0}^{T}\max_{x\in \Omega_{2}(t)}u^{4},\ea\enn owes to \eqref{2.8} and \eqref{2.4}.

Next, substituting \eqref{ll} and  \eqref{ll1} into  \eqref{as2} and utilizing the Cauchy-Schwarz inequality give rise  to
  \bnn\ba
 I_{3}
&\le C+\frac{\kappa}{8}\int_{0}^{T}\int_{\Omega} \frac{r^{2(n-1)}|\p_{x}(\theta-2)_{+}|^{2}}{v} -3\beta\int_{0}^{T}\int_{\Omega_{2}(t)} \frac{(r^{n-1}u)_{x}^{2}}{v} \\
&\quad+C\int_{0}^{T}\left[\max_{x\in \Omega}(\theta-2)_{+}^{2}+\max_{x\in \Omega_{2}(t)}\left(u^{4}+(\theta-1)^{2}\right)\right],\ea\enn
where the  $\int_{\Omega}(u^{2}+(v-1)^{2})\le C$ has been used due to  \eqref{2.4} and \eqref{e8}.

By the  Young inequality and \eqref{k5}, it satisfies
\bnn\ba I_{4}
&\le \ep \int_{0}^{T}\int_{\Omega}\th  (r^{n-1}u)_{x}^{2}+C(\ep)\int_{0}^{T}\int_{\Omega_{2}(t)} \th(\theta-2)_{+}^{2}\\
&\le \ep \int_{0}^{T}\int_{\Omega}\th  (r^{n-1}u)_{x}^{2}+C(\ep)\int_{0}^{T}\max_{x\in \Omega_{2}(t)}(\theta-2)_{+}^{2}.\ea\enn

Finally,  direct calculation from  $\eqref{eq}_{3}$  shows
  \be\ba\label{a3} I_{5}&=-2\int_{0}^{T}\int_{\Omega}  u_{t}u(\theta-2)_{+}\\
&\leq \int_{\Omega}  u_{0}^{2}(\theta_{0}-2)_{+} + \int_{0}^{T}\int_{\Omega_{2}(t)}  u^{2}\p_{t}\th\\
&\le C+\frac{\ka}{c_{v}}\int_{0}^{T}\int_{\Omega_{2}(t)}  u^{2}\left(\frac{r^{2(n-1)}\th_{x}}{v}\right)_{x}+c_{v}^{-1}\int_{0}^{T}\int_{\Omega_{2}(t)} u^{2}\tilde{R},
\ea\ee
with \bnn\ba  \tilde{R} = \beta\frac{(r^{n-1}u)_{x}^{2}}{v}- \frac{R(\theta-2)_{+}}{v}(r^{n-1}u)_{x}
-\frac{2R}{v}(r^{n-1}u)_{x}-2\mu(n-1)(r^{n-2}u^{2})_{x}.\ea\enn
In order to deal with $\int_{0}^{T}\int_{\Omega_{2}(t)}  u^{2}\left(\frac{r^{2(n-1)}\th_{x}}{v}\right)_{x},$  we define
 \bnn\ba   {\rm sgn}_{\eta}\,s=\left\{
                       \begin{array}{ll}
                        1, & s>\eta,\\
                        s/\eta, & 0\le s\le \eta,\\
                         0, &  s\le 0,
                       \end{array}
                     \right.
\ea\enn and compute by means of   Lebesgue Dominated Convergence Theorem
 \be\la{k6}\ba & \int_{0}^{T}\int_{\Omega_{2}(t)}  u^{2}\left(\frac{r^{2(n-1)}\th_{x}}{v}\right)_{x}\\
& =\lim_{\eta\rightarrow0+}\int_{0}^{T}\int_{\Omega} u^{2}{\rm sgn}_{\eta}\,(\theta-2) \left(\frac{r^{2(n-1)}\theta_{x}}{v}\right)_{x} \\
&= -\lim_{\eta\rightarrow0+}\int_{0}^{T}\int_{\Omega}\left[2uu_{x}{\rm sgn}_{\eta}\,(\theta-2)+u^{2} {\rm sgn}_{\eta}'\,(\theta-2)\right]\frac{r^{2(n-1)}\theta_{x}}{v}  \\
&\le -\lim_{\eta\rightarrow0+}\int_{0}^{T}\int_{\Omega} 2uu_{x}{\rm sgn}_{\eta}\,(\theta-2) \frac{r^{2(n-1)}\theta_{x}}{v}  \\
&\le \frac{c_{v}}{8}\int_{0}^{T}\int_{\Omega_{2}(t)}\frac{r^{2(n-1)} \th_{x}^{2}}{v}+C\int_{0}^{T}\int_{\Omega_{2}(t)}r^{2(n-1)}u^{2}u_{x}^{2}\\
&\le \frac{c_{v}}{8}\int_{0}^{T}\int_{\Omega_{2}(t)}\frac{r^{2(n-1)} \th_{x}^{2}}{v}+C\int_{0}^{T}\int_{\Omega_{2}(t)}u^{2}(r^{n-1}u)_{x}^{2}+C\int_{0}^{T}\max_{x\in \Omega_{2}(t)}u^{4},
\ea\ee where   the last inequality  owes to  \eqref{2.8},   \eqref{e8} and the simple fact
\be\la{a4} r^{n-1}u_{x}=(r^{n-1}u)_{x}-(n-1)r^{-1}vu.\ee

It is easy to check from  \eqref{2.3} that
 \be\la{1} (r^{n-2}u^{2})_{x}=ur^{-1}\left[2(r^{n-1}u)_{x}-nr^{-1}vu\right].
\ee which combining   with \eqref{z12},  \eqref{2.8} and \eqref{e8}  leads to  \bnn\ba  \int_{0}^{T}\int_{\Omega_{2}(t)}u^{2}(r^{n-2}u^{2})_{x}
 &\le C\int_{0}^{T}\int_{\Omega}u^{2}(r^{n-1}u)_{x}^{2}+ \int_{0}^{T}\max_{x\in \Omega_{2}(t)}u^{4},
\ea\enn By this, we estimate
 \be\la{k7}\ba  \int_{0}^{T}\int_{\Omega_{2}(t)}u^{2}\tilde{R}
&\le C\int_{0}^{T}\int_{\Omega}u^{2}(r^{n-1}u)_{x}^{2}+ \beta c_{v}\int_{0}^{T}\int_{\Omega } \frac{(r^{n-1}u)_{x}^{2}}{v}\\
&\quad+ C\int_{0}^{T}\left[\max_{x\in \Omega}(\th-2)_{+}^{2}+\max_{x\in \Omega_{2}(t)}u^{4}\right].
\ea\ee
With the aid of  \eqref{k6} and \eqref{k7},  the \eqref{a3} satisfies   \bnn\ba   I_{5}&\leq  C +\frac{\kappa}{8}\int_{0}^{T}\int_{\Omega} \frac{r^{2(n-1)}|\p_{x}(\theta-2)_{+}|^{2}}{v}
+C\int_{0}^{T}\int_{\Omega} u^{2}(r^{n-1}u)_{x}^{2} \\
&\quad +C\int_{0}^{T}\max_{x\in \Omega }\left(u^{4} +(\theta-2)_{+}^{2}\right)+ \beta\int_{0}^{T}\int_{\Omega}\frac{(r^{n-1}u)_{x}^{2}}{v}.
\ea\enn

On account of  estimates on  $I_{i}\,(i=1\sim5)$ above,
it follows  from \eqref{e8} and \eqref{1.3-1} that
    \be\ba\label{1.4-h}   &\int_{\Omega} (\theta-2)_{+}^{2}+   \int_{0}^{T}\int_{\Omega} \left[(r^{n-1}u)_{x}^{2}+(\th-2)_{+}(r^{n-1}u)_{x}^{2}+ r^{2(n-1)}|\p_{x}(\theta-2)_{+}|^{2} \right]\\
    &\leq  C+C\int_{0}^{T}\int_{\Omega}u^{2}(r^{n-1}u)_{x}^{2} +\ep \int_{0}^{T}\int_{\Omega}\th  (r^{n-1}u)_{x}^{2}\\
    &\quad+C\int_{0}^{T}\left(\max_{x\in \Omega}\left[(\theta-2)_{+}^{2}+ u^{4}\right]+\max_{x\in \Omega_{2}(t)}(\theta-1)^{2}\right).
\ea\ee
Noting from  \eqref{2.4}, \eqref{e8} and \eqref{a4} that
\bnn\ba   \int_{0}^{T}\int_{\Omega} r^{2(n-1)}\theta_{x}^{2}
&\le C\int_{0}^{T}\int_{\Omega} r^{2(n-1)}|\p_{x}(\theta-2)_{+}|^{2} +C\int_{0}^{T}\int_{\Omega\backslash\Omega_{2}(t)} \frac{r^{2(n-1)}\theta_{x}^{2}}{v\theta^{2}}\\
&\le C\int_{0}^{T}\int_{\Omega} r^{2(n-1)}|\p_{x}(\theta-2)_{+}|^{2}+C\ea\enn
and that \bnn\ba   \int_{0}^{T}\int_{\Omega}\th (r^{n-1} u)_{x}^{2}
&\le 3\int_{0}^{T} \int_{\Omega_{3}(t)}(\th-2)_{+} (r^{n-1} u)_{x}^{2}+3\int_{0}^{T}\int_{\Omega\backslash\Omega_{3}(t)} \frac{(r^{n-1} u)_{x}^{2}}{v\th}\\
&\le 3\int_{0}^{T} \int_{\Omega_{3}(t)}(\th-2)_{+} (r^{n-1} u)_{x}^{2}+C,\ea\enn
we select  $\ep$ in  \eqref{1.4-h} so  small such that \be\ba\label{1.4-1}   &\int_{\Omega}(\theta-2)_{+}^{2}+   \int_{0}^{T}\int_{\Omega} \left[(1+\th)(r^{n-1}u)_{x}^{2} + r^{2(n-1)}\th_{x}^{2} \right]\\
    &\leq  C+C\int_{0}^{T}\int_{\Omega}u^{2}(r^{n-1}u)_{x}^{2} \\
    &\quad+C\int_{0}^{T}\left(\max_{x\in \Omega}\left[(\theta-2)_{+}^{2}+ u^{4}\right]+\max_{x\in \Omega_{2}(t)}(\theta-1)^{2}\right).
\ea\ee

 \emph{Step 2.}
 Since \be\la{o} (r^{n-1}u^{3})_{x}=3u^{2}(r^{n-1}u)_{x}-\frac{2(n-1)}{r}vu^{3},\ee
we multiply  $\eqref{eq}_{2}$ by $u^{3}$  to receive
     \be\ba\la{lia}  &\frac{1}{4}\int_{\Omega}u^{4}(x,T)+3\beta\int_{0}^{T}\int_{\Omega} \frac{u^{2}(r^{n-1}u)_{x}^{2}}{v} \\
&=\frac{1}{4}\int_{\Omega}  u^{4}_{0}+ 2\beta(n-1) \int_{0}^{T}\int_{\Omega} \frac{(r^{n-1}u)_{x}u^{3}}{r}-\int_{0}^{T}\int_{\Omega} \left(\frac{R\theta}{v}\right)_{x}r^{n-1}u^{3}.\ea\ee

The Cauchy-Schwarz inequality, \eqref{2.4} and \eqref{2.8} show
\be\ba\la{rf}
\int_{0}^{T}\int_{\Omega} \frac{(r^{n-1}u)_{x}u^{3}}{r}\le \ep\int_{0}^{T}\int_{\Omega} (r^{n-1}u)_{x}^{2}+C(\ep)\int_{0}^{T}\max_{x\in \Omega}u^{4}(\cdot,t).\ea\ee

Next to control the last term in \eqref{lia}.  In view of   \eqref{2.4}, \eqref{e8}, and \eqref{z12},  
 \be\ba\la{o1}
&\int_{0}^{T}\int_{\Omega} \left[3u^{2}(r^{n-1}u)_{x}\frac{\th-1}{v}-2(n-1)\frac{\th-1}{r}u^{3}\right]\\
&=\int_{0}^{T}\left(\int_{\Omega\backslash \Omega_{2}(t)}+\int_{\Omega_{2}(t)} \right)\left[3u^{2}(r^{n-1}u)_{x}\frac{\th-1}{v}-2(n-1)\frac{\th-1}{r}u^{3}\right]\\
&\le \ep\int_{0}^{T}\int_{\Omega\backslash \Omega_{2}(t)} \left[(r^{n-1}u)_{x}^{2}+(\th-1)^{2}  u^{4}\right]+C\int_{0}^{T}\int_{\Omega\backslash \Omega_{2}(t)}\left[(\th-1)^{2}u^{4}+\frac{v u^{2}}{r^{2}\th}\right]\\
&\quad+\ep\int_{0}^{T}\int_{\Omega_{2}(t)}\left[u^{2}(r^{n-1}u)_{x}^{2}+(\th-1)^{2}\right]+C\int_{0}^{T}\int_{\Omega_{2}(t)}\left[(\th-1)^{2}u^{2}+u^{6}\right]\\
&\le \ep\int_{0}^{T}\int_{\Omega}(1+u^{2}) (r^{n-1}u)_{x}^{2} +C\int_{0}^{T}\left[\max_{\Omega_{2}(t)}(\th-1)^{2}+\max_{x\in \Omega}u^{4}\right]+C.
\ea\ee
Similarly, by \eqref{2.4}, \eqref{e8} and  \eqref{k5},
 \be\ba\la{o2}
&\int_{0}^{T}\int_{\Omega}\left[3u^{2}(r^{n-1}u)_{x}\frac{1-v}{v}-2(n-1)\frac{1-v}{r}u^{3}\right]\\
&\le \ep \int_{0}^{T}\int_{\Omega}(r^{n-1}u)_{x}^{2}+C \int_{0}^{T}\int_{\Omega}(1-v)^{2}u^{4}\\
&\quad+C\int_{0}^{T}\left(\int_{\Omega\backslash \Omega_{2}(t)}+\int_{ \Omega_{2}(t)}\right)\frac{|1-v|u^{3}}{r}\\
&\le \ep \int_{0}^{T}\int_{\Omega}(r^{n-1}u)_{x}^{2}+C\int_{0}^{T}\int_{\Omega} (v-1)^{2}  u^{4} \\
&\quad+C\int_{0}^{T}\int_{\Omega_{2}(t)} \th (1-v)^{2} u^{4}+C\int_{0}^{T}\int_{\Omega} \frac{v u^{2}}{r^{2}\th}\\
&\le \ep\int_{0}^{T}\int_{\Omega} (r^{n-1}u)_{x}^{2} +C\int_{0}^{T} \max_{x\in \Omega}u^{4} +C.
\ea\ee  
With the aid of  \eqref{o1} and \eqref{o2}, we infer \be\ba\la{rf1}
&-\int_{0}^{T}\int_{\Omega} \left(\frac{\theta}{v}\right)_{x}(r^{n-1}u^{3})\\
&=\int_{0}^{T}\int_{\Omega} \left(\frac{\th-1}{v} +\frac{1-v}{v} \right)(r^{n-1}u^{3})_{x}\\
&=\int_{0}^{T}\int_{\Omega}\left[3u^{2}(r^{n-1}u)_{x}\frac{\th-1}{v}-2(n-1)\frac{\th-1}{r}u^{3}\right]\\
&\quad +\int_{0}^{T}\int_{\Omega}\left[3u^{2}(r^{n-1}u)_{x}\frac{1-v}{v}-2(n-1)\frac{1-v}{r}u^{3}\right]\\
  &\le \ep\int_{0}^{T}\int_{\Omega}(1+u^{2}) (r^{n-1}u)_{x}^{2} +C\int_{0}^{T}\left[\max_{\Omega_{2}(t)}(\th-1)^{2}+\max_{x\in \Omega}u^{4}\right]+C.
\ea\ee
Choosing  $\ep$ properly small, it yields from  \eqref{rf},  \eqref{rf1}  and   \eqref{lia} that
\be\ba\la{1.5-1}
 &\int_{\Omega} u^{4}+\int_{0}^{T}\int_{\Omega} u^{2}(r^{n-1}u)_{x}^{2}\\
  &\le C+C\ep\int_{0}^{T}\int_{\Omega} (r^{n-1}u)_{x}^{2}+C\int_{0}^{T}\left(\max_{x\in \Omega} u^{4}+\max_{x\in \Omega_{2}(t)}(\th-1)^{2}\right).
    \ea\ee
Multiplying  \eqref{1.5-1} by a large constant,   adding it up to  \eqref{1.4-1}, and selecting $\ep$ small once more,   we conclude
\be\ba\label{lia5}
&\int_{\Omega} \left[(\theta-2)_{+}^{2}+ u^{4}\right]+ \int_{0}^{T}\int_{\Omega} \left[(1+\th+u^{2})(r^{n-1}u)_{x}^{2}+ r^{2(n-1)}\theta_{x}^{2}\right]\\
&\leq  C +C\int_{0}^{T}\left[\max_{x\in \Omega}(\theta-2)_{+}^{2} +\max_{x\in \Omega_{2}(t)}(\theta-1)^{2}+ \max_{x\in \Omega}u^{4}\right]\\
&\leq  C +C\int_{0}^{T}\max_{x\in \Omega}\left[(\theta-3/2)_{+}^{2} + u^{4}\right].
\ea\ee

\emph{Step 3.}
It only left to estimate  the terms in  \eqref{lia5}.

From  \eqref{k5} and \eqref{2.8}  we compute
\bnn\ba   &(\theta(x,t)-3/2)_{+}^{2}
\le C\int_{\Omega_{3/2}(t)}(\theta-3/2)_{+}|\theta_{x}|\\
&\le \frac{C}{\sqrt{\delta_{1}}} \int_{\Omega}\frac{\theta_{x}^{2}}{\theta}+\sqrt{\delta_{1}}\int_{\Omega_{3/2}(t)}(\theta-3/2)_{+}^{2} \theta \\
&\le \sqrt{\delta_{1}}\int_{\Omega}r^{2(n-1)}\theta_{x}^{2}+\frac{C}{\delta_{1}^{3/2}}\int_{\Omega} \frac{r^{2(n-1)}\theta_{x}^{2}}{v\theta^{2}}+C\sqrt{\delta_{1}}\max_{x\in \Omega} (\theta-3/2)_{+}^{2},\ea\enn
which satisfies, for  $\delta_{1}$  small, 
\bnn \max_{x\in\Omega} (\theta-3/2)_{+}^{2}\le 2\sqrt{\delta_{1}} \int_{\Omega}r^{2(n-1)}\theta_{x}^{2}
+C(\delta_{1})\int_{\Omega} \frac{r^{2(n-1)}\theta_{x}^{2}}{v\theta^{2}}.\enn
Again using \eqref{2.4},  integration it in time yields
\be\ba\la{lia9}  \int_{0}^{T}\max_{x\in \Omega}(\theta-3/2)_{+}^{2}
&\le 2\sqrt{\delta_{1}}\int_{0}^{T}\int_{\Omega}r^{2(n-1)}\theta_{x}^{2}+C(\delta_{1}).\ea\ee

Next, by \eqref{e8}, \eqref{2.8} and  \eqref{k5}, one has \bnn\ba  u^{4}(x,t)
&\le  \frac{C}{\delta_{2}}\int_{\Omega} \frac{r^{2(n-1)}u^{2}_{x}}{v\theta}+\delta_{2} \int_{\Omega}\frac{u^{6}\theta}{r^{2(n-1)}}\\
&\le \frac{C}{\delta_{2}}\int_{\Omega} \frac{r^{2(n-1)}u^{2}_{x}}{v\theta}+\delta_{2}\int_{\Omega_{2}(t)}u^{6}\theta +\delta_{2}\int_{\Omega\backslash \Omega_{2}(t)}u^{6}\theta\\
&\le \frac{C}{\delta_{2}}\int_{\Omega}\frac{r^{2(n-1)}u^{2}_{x}}{v\theta}+C\delta_{2}\max_{x\in\Omega} u^{6}(\cdot,t)+C\delta_{2}\max_{x\in\Omega} u^{4}(\cdot,t),
\ea\enn  which implies 
\be \la{ll4} \max_{x\in\Omega} u^{4}(\cdot,t)\le C(\delta_{2})\int_{\om} \frac{r^{2(n-1)}u^{2}_{x}}{v\theta}+ C\delta_{2} \max_{x\in\Omega} u^{6}(\cdot,t).
\ee
On the other hand,  from \eqref{a4} and \eqref{2.4} we compute \bnn\ba  \max_{x\in \Omega} u^{6}(\cdot,t)
& = 6\int_{0}^{x}\frac{u^{5} [(r^{n-1}u)_{x}-(n-1) r^{-1}vu]}{r^{n-1}}\\
&\le  6\int_{0}^{x}\frac{u^{5} (r^{n-1}u)_{x}}{r^{n-1}}\\
&\le  \frac{C}{\sqrt{\delta_{2}}}\int_{\Omega}u^{2}(r^{n-1}u)_{x}^{2} +C\sqrt{\delta_{2}}\max_{x\in\Omega} u^{6}(\cdot,t).
\ea\enn
This,  together with \eqref{ll4} and \eqref{2.4},  yields
\be \ba\la{lia21}  \int_{0}^{T}\max_{x\in \Omega} u^{4}(\cdot,t)\le C(\delta_{2}) + C\sqrt{\delta_{2}} \int_{0}^{T}\int_{\Omega}u^{2}(r^{n-1}u)_{x}^{2}.\ea\ee

Substituting  \eqref{lia21} and \eqref{lia9}    back into \eqref{lia5} gives birth to 
  \bnn\ba
 \int_{\Omega}  \left[(\theta-2)_{+}^{2}+ u^{4}\right]+\int_{0}^{T} \int_{\Omega}\left[(1+\th+ u^{2})(r^{n-1}u)_{x}^{2} + r^{2(n-1)}\theta_{x}^{2}\right]
    &\leq  C, \ea\enn as long as $\delta_{1}$ and $\delta_{2}$ are chosen small enough.  
    
 Finally, 
\bnn\ba \int_{\Omega} (\theta-1)^{2}= \left(\int_{\Omega \backslash\Omega_{3}(t)}+\int_{ \Omega_{3}(t)}\right)  (\theta-1)^{2}
 \le C+C\int_{\Omega_{3}(t)} (\theta-2)_{+}^{2},\ea\enn by \eqref{2.4}.
In conclusion, the  last two inequalities  guarantee  \eqref{e}, the required. $\Box$

 \begin{lemma} \label{l5.1}  For the solution $(v,u,\th)$  described in Theorem \ref{t1}, it holds
\be\la{e1}  \sup_{0\le t\le T}\int_{\Omega}v_{x}^{2}(x,t)+ \int_{0}^{T} \int_{\Omega} (1+\th)v_{x}^{2}
  \leq C.\ee
 \end{lemma}
\emph{Proof.}  By $\eqref{eq}_{1}$, rewriting   $\eqref{eq}_{2}$ as the form
\bnn \beta\left(\frac{ v_{x}}{v}\right)_{t}=R \left(\frac{\th}{v}\right)_{x}+ r^{1-n}u_{t},\enn
which yields after multiplied  by $v_{x}/v$
\be\ba\la{hy} &\frac{\beta}{2}\int_{\Omega}\frac{ v_{x}^{2}}{v^{2}}(x,T) +\int_{0}^{T}\int_{\Omega}\frac{R\th v_{x}^{2}}{v^{3}}\\
& = \frac{\beta}{2}\int_{\Omega}\frac{ v_{x}^{2}}{v^{2}}(x,0)+ \int_{0}^{T}\int_{\Omega}\frac{Rv_{x}\th_{x}}{v^{2}}+\int_{0}^{T}\int_{\Omega} r^{1-n}u_{t}\frac{ v_{x}}{v}.\ea\ee

The Cauchy-Schwarz inequality, \eqref{2.4}, \eqref{2.8}  and \eqref{e} ensure that
  \bnn\ba  \int_{0}^{T}\int_{\Omega}\frac{Rv_{x}\th_{x}}{v^{2}}&\le \frac{1}{2} \int_{0}^{T}\int_{\Omega}\frac{R\th  v_{x}^{2}}{v^{3}}+C  \int_{0}^{T}\int_{\Omega}\frac{\th_{x}^{2}}{v\th}\\
&\le \frac{1}{2} \int_{0}^{T}\int_{\Omega}\frac{R\th  v_{x}^{2}}{v^{3}}+C\int_{0}^{T}\int_{\Omega}\frac{r^{2(n-1)}\th_{x}^{2}}{v\th^{2}}+C\int_{0}^{T}\int_{\Omega}r^{2(n-1)}\th_{x}^{2}\\
&\le  \frac{1}{2} \int_{0}^{T}\int_{\Omega}\frac{R\th  v_{x}^{2}}{v^{3}}+C.\ea\enn
Thanks to  \eqref{2.3}, \eqref{e2}, and \eqref{2.4},
 \bnn\ba &\int_{0}^{T}\int_{\Omega} r^{1-n}u_{t}\frac{ v_{x}}{v}\\
  &=\int_{\Omega} r^{1-n}u \frac{ v_{x}}{v}(x,T)-\int_{\Omega} r^{1-n}u \frac{ v_{x}}{v}(x,0)\\
  &\quad-\int_{0}^{T}\int_{\Omega} r^{1-n}u\left(\frac{ v_{x}}{v}\right)_{t}+(n-1)\int_{0}^{T}\int_{\Omega} r^{-n}u^{2}\frac{ v_{x}}{v} \\
  &\le C+\frac{\beta}{4}\int_{\Omega}\left|\frac{ v_{x}}{v} \right|^{2}-\int_{0}^{T}\int_{\Omega} r^{1-n}u\left(\frac{ v_{x}}{v}\right)_{t}+(n-1)\int_{0}^{T}\int_{\Omega} r^{-n}u^{2}\frac{ v_{x}}{v}.\ea\enn
Hence, the  \eqref{hy} satisfies 
 \be\ba\la{hyu} &\frac{\beta}{4}\int_{\Omega}\left|\frac{ v_{x}}{v} \right|^{2}+\frac{1}{2} \int_{0}^{T}\int_{\Omega}\frac{R\th v_{x}^{2}}{v^{3}}\\
& \le C-\int_{0}^{T}\int_{\Omega} r^{1-n}u\left(\frac{ v_{x}}{v}\right)_{t}+(n-1)\int_{0}^{T}\int_{\Omega} r^{-n}u^{2}\frac{ v_{x}}{v}.\ea\ee

In terms of $\eqref{eq}_{1}$, \eqref{2.8},  \eqref{e}, \eqref{e8}, as well as 
 \bnn\ba (r^{1-n}u)_{x}=r^{2(1-n)}(r^{n-1}u)_{x}+2(1-n)r^{1-2n} vu,\ea\enn
it satisfies  \be\la{k}\ba  &-\int_{0}^{T}\int_{\Omega} r^{1-n}u\left(\frac{ v_{x}}{v}\right)_{t}\\
&=\int_{0}^{T}\int_{\Omega}(r^{1-n}u)_{x} \frac{(r^{n-1}u)_{x}}{v}\\
&=\int_{0}^{T}\int_{\Omega} r^{2(1-n)}\frac{(r^{n-1}u)_{x}^{2}}{v}+2(1-n)\int_{0}^{T}\int_{\Omega} r^{1-2n} u(r^{n-1}u)_{x}\\
&\le  C\int_{0}^{T}\int_{\Omega} (1+\th)(r^{n-1}u)_{x}^{2} +C\int_{0}^{T}\int_{\Omega}\frac{v u^{2}}{r^{2}\th}\le C.\ea\ee
From  \eqref{e8} and \eqref{3.12} we have\be\la{xx}\int_{0}^{T}\int_{\Omega}\frac{v_{x}^{2}}{v^{2}}
 \le C(\alpha_{1},\alpha_{2})\left(\int_{0}^{T}f(t)\int_{\Omega}\frac{v_{x}^{2}}{v^{2}} + \int_{0}^{T}\int_{\Omega}\frac{\th v_{x}^{2}}{v^{2}}\right), \ee
where  $f(t)$ is defined in  \eqref{a6}.  Therefore,   \be\la{k1}\ba  \int_{0}^{T}\int_{\Omega} r^{-n}u^{2}\frac{ v_{x}}{v}
  &\le C \int_{0}^{T}\max_{x\in \om} u^{4}\int_{\Omega} r^{-2n}+ \frac{1}{4}\left(\int_{0}^{T}f(t)\int_{\Omega}\frac{v_{x}^{2}}{v^{2}}+ \int_{0}^{T}\int_{\Omega}\frac{R\th v_{x}^{2}}{v^{2}}\right)\\
   &\le C + \frac{1}{4}\left(\int_{0}^{T}f(t)\int_{\Omega}\frac{v_{x}^{2}}{v^{2}}+ \int_{0}^{T}\int_{\Omega}\frac{R\th v_{x}^{2}}{v^{2}}\right),\ea\ee
where we also used  \eqref{o4}, \eqref{lia21} and \eqref{e}.  Substituting   \eqref{k1} and \eqref{k}  into  \eqref{hyu}   arrives at
 \bnn\ba  \int_{\Omega}\left|\frac{ v_{x}}{v} \right|^{2}+ \int_{0}^{T}\int_{\Omega}\frac{\th v_{x}^{2}}{v^{3}} \le C +C\int_{0}^{T}f(t)\int_{\Omega}\frac{v_{x}^{2}}{v^{2}}.\ea\enn This and \eqref{xx}   complete the proof after integration in time.  $\Box$

\begin{lemma} \label{l5.2} It holds that
\be\la{f}  \sup_{0\le t\le T}\int_{\Omega}u_{x}^{2}(x,t)+ \int_{0}^{T} \int_{\Omega} r^{2(n-1)}u_{xx}^{2}
  \leq C\left(1+\max_{\Omega\times [0,T]}\th \right).\ee  \end{lemma}
\emph{Proof.} Multiplied by $-u_{xx}$,  it yields from $\eqref{eq}_{2}$  that
\bnn\ba & \frac{1}{2}\p_{t} u_{x}^{2} +\beta  \frac{r^{2(n-1)}u_{xx}^{2}}{v}\\
&=(u_{x}u_{t})_{x}+\beta u_{xx}\left(r^{2(n-1)}\frac{v_{x}u_{x}}{v^{2}}+(n-1)\frac{uv}{r^{2}}-2(n-1)r^{n-2}u_{x}\right)\\
&\quad+Ru_{xx}r^{n-1}\left(\frac{\th_{x}}{v}- \frac{\th v_{x}}{v^{2}}\right).\ea\enn
Integrating the above equality leads  to
\be\la{o5}\ba & \frac{1}{2}\int_{\Omega}u_{x}^{2}(x,T)+\beta \int_{0}^{T}\int_{\Omega} \frac{r^{2(n-1)}u_{xx}^{2}}{v}\\
&\le \frac{1}{2}\int_{\Omega}u_{0x}^{2}+\frac{\beta}{4}\int_{0}^{T}\int_{\Omega} \frac{r^{2(n-1)}u_{xx}^{2}}{v}\\
&\quad +C\int_{0}^{T}\int_{\Omega}\left[ r^{2(n-1)} v_{x}^{2}u_{x}^{2} +\frac{u^{2}}{r^{2(n+1)}} + \frac{u_{x}^{2}}{r^{2}}+\th_{x}^{2}+\th^{2}v_{x}^{2}\right].\ea\ee
By virtue of  \eqref{e8}, \eqref{e},   \eqref{e1},   \eqref{2.8} and  \eqref{2.4},
\bnn\ba &C\int_{0}^{T}\int_{\Omega}\left(\frac{u^{2}}{r^{2(n+1)}} + \frac{u_{x}^{2}}{r^{2}}+\th_{x}^{2}+\th^{2}v_{x}^{2}\right)\\
&\le C\left(1+\max_{\Omega\times [0,T]}\th \right) \int_{0}^{T}\int_{\Omega}\left(\frac{v u^{2}}{r^{2}\th} +  \frac{r^{2(n-1)}u_{x}^{2}}{v\th}+r^{2(n-1)}\th_{x}^{2}+\th v_{x}^{2}\right)\\
&\le C\left(1+\max_{\Omega\times [0,T]}\th \right).
 \ea\enn
Since $H^{1}\hookrightarrow L^{\infty}$, we use \eqref{e1} and  \eqref{2.4} to get
\be\la{o7}\ba  C\int_{0}^{T}\int_{\Omega} r^{2(n-1)} v_{x}^{2}u_{x}^{2}
&\le C \int_{0}^{T}\|r^{n-1}u_{x}\|_{L^{\infty}}^{2} \\
&\le  \frac{\beta}{4}\int_{0}^{T}\int_{\Omega}  r^{2(n-1)}u_{xx}^{2}+C \int_{0}^{T}\int_{\Omega} r^{2(n-1)}u_{x}^{2} \\
&\le  \frac{\beta}{4}\int_{0}^{T}\int_{\Omega}  r^{2(n-1)}u_{xx}^{2}+C \max_{\Omega\times [0,T]}\th.
 \ea\ee
 With the last three inequalities in hand,  we get the \eqref{f}.    $\Box$

\begin{lemma} \label{l5.3} It holds that
\be\la{f1} \sup_{0\le t\le T} \int_{\Omega}\th_{x}^{2}(x,t)+\int_{0}^{T} \int_{\Omega} r^{2(n-1)}\th_{xx}^{2}
  \leq  C\left(1+\max_{\Omega\times [0,T]}\th^{2} \right).\ee
 \end{lemma}
\emph{Proof.} Multiplying $\eqref{eq}_{3}$ by $-\th_{xx}$ yields
\bnn\ba &\frac{c_{v}}{2}\p_{t} \th_{x}^{2} +\ka  \frac{r^{2(n-1)}\th_{xx}^{2}}{v}\\
&=(c_{v}\th_{x}\th_{t})_{x}+\ka \th_{xx}\left(r^{2(n-1)}\frac{v_{x}\th_{x}}{v^{2}} -2(n-1)r^{n-2}\th_{x}\right)\\
&\quad+\th_{xx}\left(\frac{R\th}{v}(r^{n-1}u)_{x}+2\mu(n-1)(r^{n-2}u^{2})_{x}-\beta \frac{(r^{n-1}u)_{x}^{2}}{v}\right).\ea\enn
Utilizing  the  Cauchy-Schwarz inequality,   \eqref{a4} and \eqref{1}, we integrate it in time and deduce 
\be\la{o8}\ba & \frac{c_{v}}{2}\int_{\Omega}\th_{x}^{2}+\ka\int_{0}^{T}\int_{\Omega} \frac{r^{2(n-1)}\th_{xx}^{2}}{v}\\
&\le \frac{c_{v}}{2}\int_{\Omega}\th_{0x}^{2}+\frac{\ka}{4}\int_{0}^{T}\int_{\Omega}\frac{r^{2(n-1)}\th_{xx}^{2}}{v}
+C\int_{0}^{T}\int_{\Omega}  \frac{r^{2(n-1)}v_{x}^{2}\th_{x}^{2}}{v^{3}}\\
&\quad +C\int_{0}^{T}\int_{\Omega}\left[\th_{x}^{2} + (\th^{2}+r^{-2}u^{2})(r^{n-1}u)_{x}^{2}+r^{-2(n+1)}u^{4}+r^{2(n-1)}u_{x}^{4} \right].\ea\ee
Thanks  to  \eqref{e1} and    \eqref{e},  a similar argument as \eqref{o7} shows
\bnn\ba  C\int_{0}^{T}\int_{\Omega}  \frac{r^{2(n-1)}v_{x}^{2}\th_{x}^{2}}{v^{3}}
&\le  \frac{\ka}{4}\int_{0}^{T}\int_{\Omega}  \frac{r^{2(n-1)}\th_{xx}^{2}}{v}+C\int_{0}^{T}\int_{\Omega} r^{2(n-1)}\th_{x}^{2} \\
&\le \frac{\ka}{4}\int_{0}^{T}\int_{\Omega} \frac{r^{2(n-1)}\th_{xx}^{2}}{v}+C.
 \ea\enn
 Inequalities \eqref{2.8}, \eqref{o4}, \eqref{lia21}, and  \eqref{e} guarantee that
\bnn\ba  &\int_{0}^{T}\int_{\Omega}\left( \th_{x}^{2} + (\th^{2}+r^{-2}u^{2})(r^{n-1}u)_{x}^{2}+r^{-2(n+1)}u^{4} \right) \\
&\le  C\left(1+\max_{\Omega\times [0,T]}\th \right)\int_{0}^{T}\int_{\Omega}\left[r^{2(n-1)}\th_{x}^{2} + (\th+u^{2}) (r^{n-1}u)_{x}^{2} \right]\\
&\quad+C\int_{0}^{T}\max_{x\in \Omega}u^{4}\int_{\Omega}r^{-2(n+1)} \\
&\le  C\left(1+\max_{\Omega\times [0,T]}\th \right).\ea\enn
Finally,  from    \eqref{f} and  \eqref{o7} we obtain
\bnn\ba \int_{0}^{T}\int_{\Omega}  r^{2(n-1)}u_{x}^{4}
&\le C\int_{0}^{T}\|r^{n-1}u_{x}\|_{L^{\infty}}^{2}\int_{\Omega} u_{x}^{2} \\
&\le C\max_{\Omega\times [0,T]}\th\int_{0}^{T}\|r^{n-1}u_{x}\|_{L^{\infty}}^{2} \le C\left(1+\max_{\Omega\times [0,T]}\th^{2} \right).  \ea\enn
Inserting  the last three  inequalities into  \eqref{o8} yields  \eqref{f1}.  $\Box$

\begin{corollary}[Upper bound on  $u$ and $\th$]  For all  $(x,t)\in \overline{\Omega}\times [0,\infty)$, it satisfies \be\la{a9} |u(x,t)|+ \th(x,t) \le C.\ee
 \end{corollary}
\emph{Proof}.
In view of \eqref{e} and  \eqref{f1},   Sobolev inequality implies
\bnn \|\th-1\|_{L^{\infty}(\om)}^{2}\le C\|\th-1\|_{L^{2}(\Omega)}\|\th_{x}\|_{L^{2}(\Omega)}\le C\left(1+\max_{\Omega\times [0,T]}\th\right),\enn
which, for some $C$  independent of $T$,
\be\la{w2} \th(x,t) \le C,\quad (x,t)\in \overline{\Omega}\times [0,T].\ee   Having \eqref{w2} in hand, we conclude from \eqref{2.4} and \eqref{f} that
  \bnn  |u(x,t)| \le C,\quad (x,t)\in \overline{\Omega}\times [0,T].\enn

\section{Estimates on  derivatives}
The  rest lemmas   are  concerned with  the derivative   estimates, which are needed to show    the large-time behavior of solutions.
\begin{lemma} \label{l5.5} Let $(v,u,\th)$ be the solution obtained  in Theorem \ref{t1}.  Then
\be\la{f2}  \sup_{0\le t\le T}\int_{\Omega}r^{2(n-1)}v_{x}^{2}(x,t)+\int_{0}^{T} \int_{\Omega}(1+\th)r^{2(n-1)}v_{x}^{2}\le C.\ee
 \end{lemma}
\emph{Proof}. By \eqref{2.3} and $\eqref{eq}_{1}$, it follows  from   $\eqref{eq}_{2}$  that
\bnn \beta\left(r^{n-1}\frac{ v_{x}}{v}\right)_{t}+ r^{n-1} \frac{R\th v_{x}}{v^{2}}=   r^{n-1} \frac{R\th_{x}}{v}+u_{t} -\beta(n-1)r^{n-2} u\frac{v_{x}}{v},\enn
which gives, after multiplied by $r^{n-1}v_{x}/v,$
\be\la{a7q}\ba &\frac{\beta}{2}\int_{\Omega} \left|r^{n-1}\frac{v_{x}}{v}\right|^{2}(x,T)
+ \int_{0}^{T}\int_{\Omega} r^{2(n-1)} \frac{R\th v_{x}^{2}}{v^{3}}\\
&\le C+\frac{1}{4}\int_{0}^{T}\int_{\Omega} r^{2(n-1)} \frac{R\th v_{x}^{2}}{v^{3}}+C \int_{0}^{T}\int_{\Omega} r^{2(n-1)}\frac{\th_{x}^{2}}{v\th} \\
&\quad +\int_{0}^{T}\int_{\Omega} u_{t}r^{n-1}\frac{v_{x}}{v}-\beta (n-1)\int_{0}^{T}\int_{\Omega}r^{2(n-1)-1}u\frac{v_{x}^{2}}{v^{2}}.
\ea\ee
By virtue of   $\eqref{eq}_{1}$, \eqref{2.4}, \eqref{e} and \eqref{a9}, it satisfies
\be\la{a8}\ba  \int_{0}^{T}\int_{\Omega} u_{t}r^{n-1}\frac{v_{x}}{v}
&= \int_{\Omega} u r^{n-1}\frac{v_{x}}{v}(x,T)-\int_{\Omega} u r^{n-1}\frac{v_{x}}{v}(x,0)\\
&\quad-(n-1)\int_{0}^{T}\int_{\Omega}  u^{2} r^{n-2} \frac{v_{x}}{v}+\int_{0}^{T}\int_{\Omega} \frac{(r^{n-1}u)_{x}^{2}}{v}\\
&\le C+\frac{\beta}{4}\int_{\Omega} \left|r^{n-1}\frac{v_{x}}{v}\right|^{2} +\frac{1}{8}\int_{0}^{T}\int_{\Omega} r^{2(n-1)} \frac{R\th v_{x}^{2}}{v^{3}}.
\ea\ee
 Sobolev inequality, \eqref{2.8}, \eqref{e8}  and \eqref{a9}  ensure that 
\bnn\ba \|r^{-1}u\|_{L^{\infty}}^{2}&\le C\left(\|ur^{-1}\|_{L^{2}(\Omega)}^{2}+\|(ur^{-1})_{x}\|_{L^{2}(\Omega)}^{2}\right)  \\
&\le C\int_{\Omega}\frac{vu^{2}}{r^{2}\th}+C\int_{\Omega}\frac{r^{2(n-1)}u_{x}^{2}}{v\th}=:g(t).\ea\enn
By this,  the   Cauchy-Schwarz inequality and  \eqref{xx} show 
\be\la{w6}\ba & \int_{0}^{T}\int_{\Omega}r^{2(n-1)-1}u\frac{v_{x}^{2}}{v^{2}}\\
&\le  \int_{0}^{T}\left(\ep+C(\ep) \|r^{-1}u\|_{L^{\infty}}^{2}\right) \int_{\Omega} r^{2(n-1)}\frac{v_{x}^{2}}{v^{2}}\\
&\le  C\ep   \int_{0}^{T}\int_{\Omega}r^{2(n-1)}\frac{\th v_{x}^{2}}{v^{2}}+C(\ep)\int_{0}^{T}\left(f(t)+g(t)\right)\int_{\Omega}r^{2(n-1)}\frac{v_{x}^{2}}{v^{2}}.
\ea\ee

For  sufficiently small   $\ep,$ we
substitute  \eqref{a8} and \eqref{w6} into   \eqref{a7q}, and use \eqref{e},  \eqref{e8}, \eqref{2.4}, Gronwall's inequality, to discover

 \bnn   \int_{\Omega} r^{2(n-1)}v_{x}^{2}
+\int_{0}^{T}\int_{\Omega} r^{2(n-1)} \th v_{x}^{2} \le C.
\enn
This inequality and   \eqref{xx}   give  birth to \eqref{f2}.  $\Box$
\begin{lemma} \label{l5.6} Let $(v,u,\th)$ be the solution obtained in Theorem \ref{t1}. It holds that
\be\la{4}   \sup_{0\le t\le T}\int_{\Omega}(r^{n-1}u)_{x}^{2}(x,t)+\int_{0}^{T} \int_{\Omega}u_{t}^{2}\le C.\ee
  \end{lemma}
\emph{Proof}. Multiplying $\eqref{eq}_{2}$ by $u_{t}$ yields
\bnn\ba  \int_{0}^{T}\int_{\Omega} u_{t}^{2} &+\beta\int_{0}^{T}\int_{\Omega} \frac{(r^{n-1}u)_{x}}{v}(r^{n-1}u_{t})_{x}\\
&=\int_{0}^{T}\int_{\Omega}\left( \frac{R\th v_{x}}{v^{2}}-\frac{R\th_{x}}{v}\right) r^{n-1}u_{t} \le  \frac{1}{2} \int_{0}^{T}\int_{\Omega} u_{t}^{2}+C,
\ea\enn
owes to  \eqref{e8},  \eqref{a9}, \eqref{e}  and \eqref{f2}.

Recalling from \eqref{1}, \eqref{2.8},  \eqref{e8}, \eqref{a9}, \eqref{e}, and  \eqref{2.4},      we compute
\bnn\ba   & \int_{0}^{T}\int_{\Omega} \frac{(r^{n-1}u)_{x}}{v}(r^{n-1}u_{t})_{x} \\
&= \int_{0}^{T}\int_{\Omega} \frac{(r^{n-1}u)_{x}}{v}[(r^{n-1}u)_{t}-(n-1)(r^{n-2}u^{2})]_{x}\\
&=-\frac{1}{2}  \int_{\Omega} \frac{(r^{n-1}u)_{x}^{2}}{v}(x,0)+\frac{1}{2}  \int_{\Omega} \frac{(r^{n-1}u)_{x}^{2}}{v}(x,T)\\
&\quad+\frac{1}{2} \int_{0}^{T}\int_{\Omega} \frac{(r^{n-1}u)_{x}^{3}}{v^{2}}- (n-1)\int_{0}^{T}\int_{\Omega} \frac{(r^{n-1}u)_{x}}{v}(r^{n-2}u^{2})_{x} \\
&\ge \frac{1}{2}  \int_{\Omega} \frac{(r^{n-1}u)_{x}^{2}}{v}(x,t)-C\int_{0}^{T}\int_{\Omega} (r^{n-1}u)_{x}^{4}-C.
\ea\enn
Then,
\be\la{3}\ba \int_{0}^{T}\int_{\Omega} u_{t}^{2}+  \int_{\Omega}  (r^{n-1}u)_{x}^{2}  \le C+C\int_{0}^{T}\|(r^{n-1}u)_{x}\|_{L^{\infty}}^{2}\int_{\Omega} (r^{n-1}u)_{x}^{2} .
\ea\ee
It follows  from \eqref{2.8}, \eqref{a4},  \eqref{e1},   \eqref{f}, \eqref{2.4}  and  \eqref{a9}    that
 \be\ba\la{6}\int_{0}^{T}\|(r^{n-1}u)_{x}\|_{L^{\infty}}^{2}
 &\le C\int_{0}^{T}\left(\| (r^{n-1}u)_{xx} \|_{L^{2}(\Omega)}^{2}+ \| (r^{n-1}u)_{x}\|_{L^{2}(\Omega)}^{2}\right)\\
&\le C\int_{0}^{T} \int_{\Omega}\left(r^{2(n-1)}u_{xx}^{2} +v_{x}^{2}+\frac{r^{2(n-1)}u_{x}^{2}}{v\th}+\frac{v u^{2}}{r^{2}\th}\right)\le C.\ea\ee
By \eqref{6}, applying  Gronwall's  inequality  to \eqref{3}   receives   the  \eqref{4}. $\Box$

\begin{lemma} \label{l5.7} Let $(v,u,\th)$ be the solution obtained in Theorem \ref{t1}. Then it holds that
\be\la{7} \sup_{0\le t\le T}\int_{\Omega}r^{2(n-1)}\th_{x}^{2}(x,t)+\int_{0}^{T} \int_{\Omega}\th_{t}^{2}\le C.\ee\end{lemma}
\emph{Proof}.  Multiplied by   $\th_{t}$, it gives from  $\eqref{eq}_{3}$ that
\bnn\ba   &c_{v}\int_{0}^{T}\int_{\Omega} \th_{t}^{2}+\frac{\ka}{2}  \int_{\Omega} \frac{r^{2(n-1)}\th_{x}^{2}}{v}\\
&=\frac{\ka}{2}  \int_{\Omega} \frac{r^{2(n-1)}\th_{x}^{2}}{v}(x,0)+\frac{\ka}{2}\int_{0}^{T}\int_{\Omega}\left( 2(n-1)\frac{r^{2(n-1)}r^{-1}u\th_{x}^{2}}{v}- \frac{r^{2(n-1)}\th_{x}^{2}v_{t}}{v^{2}}\right)\\
&\quad+ \beta\int_{0}^{T}\int_{\Omega} \frac{(r^{n-1}u)_{x}^{2}\th_{t}}{v}- \int_{0}^{T}\int_{\Omega} \frac{R\th}{v}(r^{n-1}u)_{x}\th_{t}-2(n-1)\mu\int_{0}^{T}\int_{\Omega} \left(r^{n-2}u^{2}\right)_{x}\th_{t},
\ea\enn
which deduces from  the Cauchy-Schwarz inequality,   \eqref{e8}, \eqref{e}, \eqref{a9}, \eqref{2.8},  $\eqref{eq}_{1}$, \eqref{2.4},   \eqref{1}  and \eqref{4} that
\bnn  \int_{0}^{T}\int_{\Omega} \th_{t}^{2}+ \int_{\Omega} r^{2(n-1)}\th_{x}^{2}
 \le C+ C\int_{0}^{T} \|(r^{n-1}u)_{x}\|_{L^{\infty}(\Omega)}^{2}\left(1+\int_{\Omega}r^{2(n-1)}\th_{x}^{2} \right).
\enn
Hence, the \eqref{7} comes from   Gronwall's inequality and \eqref{6}.

\section{Proof of Theorem \ref{t2}}

We are now in a position to prove  Theorem \ref{t2}.

Firstly, \eqref{a2} follows from   \eqref{2.4}, \eqref{e8}, eqs \eqref{eq} and  Lemmas \ref{lem1},   \ref{l3.2},   \ref{l5.5}-\ref{l5.7}.

In view of  \eqref{2.4}, \eqref{e8}  and  \eqref{e}, to prove \eqref{a03} it suffices to justify 
 \be\la{b3}\ba  \lim_{t\rightarrow\infty}\|(v_{x}, u_{x}, \th_{x})(\cdot,t)\|_{L^{2}(\Omega)}=0.
\ea\ee
To this end, we claim
\be\la{w3}\ba
  &\int_{0}^{\infty}\left| \frac{d}{dt}\|v_{x}\|_{L^{2}(\Omega)}^{2}\right|+\left| \frac{d}{dt}\|u_{x}\|_{L^{2}(\Omega)}^{2}\right|+\left|\frac{d}{dt}\|\th_{x}\|_{L^{2}(\Omega)}^{2}\right|dt \le C.\ea\ee
In fact,  by  \eqref{2.8}, Lemmas  \ref{l5.2}-\ref{l5.3} and  Lemmas \ref{l5.6}-\ref{l5.7}, we have
 \bnn \ba &\int_{0}^{\infty} \left|\frac{d}{dt}\|u_{x}\|_{L^{2}(\Omega)}^{2}\right|+\left|\frac{d}{dt}\|\th_{x}\|_{L^{2}(\Omega)}^{2}\right|dt\\
& \le \int_{0}^{\infty} \left(\|u_{xx}\|_{L^{2}(\Omega)}^{2}+\|u_{t}\|_{L^{2}(\Omega)}^{2}+\|\th_{xx}\|_{L^{2}(\Omega)}^{2}+\|\th_{t}\|_{L^{2}(\Omega)}^{2}\right)dt\le C.\ea\enn
By $\eqref{eq}_{1}$,  \eqref{e1} and  \eqref{6},  one has
\bnn\ba \int_{0}^{\infty} \left|\frac{d}{dt}\|v_{x}\|_{L^{2}(\Omega)}^{2}\right|dt\le C\int_{0}^{\infty} \left(\|v_{x}\|_{L^{2}(\Omega)}^{2}+\|(r^{n-1}u)_{xx}\|_{L^{2}(\Omega)}^{2}\right)dt\le C.\ea\enn
 The validity of  \eqref{b3} is guaranteed by   \eqref{w3} and  the following inequality
\bnn \int_{0}^{\infty}\left(\|v_{x}\|_{L^{2}(\Omega)}^{2}+\|u_{x}\|_{L^{2}(\Omega)}^{2}+\|\th_{x}\|_{L^{2}(\Omega)}^{2}\right)dt\le C,\enn
which  comes from  \eqref{e1},   \eqref{e}, \eqref{2.4},  \eqref{e8}, \eqref{a9} and \eqref{2.8}.

In terms of \eqref{e8} and \eqref{a9},  to finish   \eqref{a1},   it remains  to show  that $\th(x,t)$  is  positive bounded from below.  Recalling \eqref{a03},   for some  $T_{1}$ large enough, one has  
\bnn \th(x,t)\ge 1/2, \quad \forall\,\,\,(x,t)\in \overline{\Omega}\times [T_{1},\infty).\enn
While for the interval $[0,T_{1}]$,   it satisfies  from \cite[eq.(4.9)]{jiang1} that
\bnn \th(x,t)\ge C(T_{1}),\quad (x,t)\in \overline{\Omega}\times [0,T_{1}].\enn
Thus, the proof of 
   Theorem \ref{t2} is complete.

\end{document}